\documentclass[12pt]{amsart}
\newtheorem{teo}{Theorem}[section]

\newtheorem{lemm}{Lemma}[section]
\newtheorem{remark}{Remark}[section]

\renewcommand{\l}{\lambda}

\newcommand{\beq}{\begin{equation}}
\newcommand{\eeq}{\end{equation}}
\newcommand{\bqn}{\begin{eqnarray}}
\newcommand{\eqn}{\end{eqnarray}}
\newcommand{\bqne}{\begin{eqnarray*}}
\newcommand{\eqne}{\end{eqnarray*}}

\newcommand{\R}{{\mathbb R}}

\newcommand{\C}{{\mathbb C}}

\renewcommand{\l}{[}
\renewcommand{\r}{]}

\begin{document}

\title[On some properties of the manifolds with skew-symmetric
torsion]{On some properties of the manifolds with skew-symmetric
torsion\\ and holonomy SU(n) and Sp(n)}

\author{Anna Fino
\and
  Gueo Grantcharov
  }

\maketitle

\noindent{{\bf Abstract} In this paper we provide examples of
hypercomplex manifolds  which do not carry HKT structure,
thus answering a question in \cite{GP}. We also prove that the
existence of HKT structure is not stable under small deformations.
Similarly we provide examples of compact complex manifolds with
vanishing first Chern class which do not admit a Hermitian
structure with restricted holonomy of its Bismut connection in
$SU(n)$, thus providing a counter-example of the conjecture in
\cite{GIP}. Again we prove that such property is not stable under
small deformations.

\section{Introduction}

Let $(M,J,g)$ be a Hermitian manifold. Then there is a unique
Hermitian connection $\nabla$ with torsion $T$ such that
$g(X,T(Y,Z))$ is totally skew-symmetric. This connection was used
by Bismut to prove a local index formula for the Dolbeault
operator when the manifold is not K\"ahler. Metric connections
with such torsion have also applications in type II string theory,
supersymmetric $\sigma$-models and the geometry of black hole
moduli spaces - see e.g. \cite{Pap} for more references. For some
applications the holonomy of the connection is restricted to a
proper subgroup of $U(n)$. We concentrate here on $SU(n)$ and
$Sp(n)$ as a possible (restricted) holonomy groups.

     In (4,0) supersymmetric $\sigma$-models with Wess-Zumino term the
     target manifold has a metric connection with skew-symmetric
     torsion with holonomy in $Sp(n)$. Such manifolds carry a triple of
anticommuting complex
     structures called hypercomplex structure compatible with the metric and
     are known in the
     physics literature as HKT manifolds \cite{HP}. A lot of information about
     their geometry is known (from a mathematical viewpoint see \cite{GP}).
      In particular it is known that there is a local existence result,
reduction theory, non-trivial Dolbeault cohomology properties and
many examples.
     So a natural question is (see \cite{GP}) {\it Do all hypercomplex
manifolds admit HKT
     structure?} Moreover for some basic examples as $SU(3)$ and
     $S^1\times S^{4n-1}$ there is also a local deformation rigidity -
      any small deformation of the invariant hypercomplex structure
      there admits HKT structure \cite{GP, GPP}. This was generalized
by M. Verbitsky \cite{V1}
     to all manifolds for which $\partial$-closed (2,0)-forms are
(2,0)-components of a closed
     forms. It rises also the question {\it Whether
      every small hypercomplex deformation of HKT structure is HKT?}
      One of the purposes of this paper is to provide negative answer
      to both questions.

      When the holonomy group of connection with skew-symmetric torsion  on
      Hermitian manifold is contained in $SU(n)$ then the manifold has
vanishing first Chern class. Such models in string theory has been
considered first by A. Strominger \cite{St} and C. Hull \cite{Hu}.
On a large class of manifolds however, the models suggested in
\cite{St} are degenerate (with vanishing torsion) \cite{IP,GMPW}.
More
    recently $SU(3)$-structures on non-K\"ahler manifolds have
attracted attention
   as a more general models for string compactifications and many examples have
been found:  \cite{BD,BBDG,CCDLMZ,GMW,GoP}. It was conjectured in \cite{GIP}
that {\it any} compact complex manifold with vanishing first Chern
class admits a Hermitian metric and connection with totally
skew-symmetric torsion and (restricted) holonomy in $SU(n)$.
      The other aim of this paper is to provide counter-example of
      this conjecture. Moreover we show that, as in the hypercomplex
      case, there is no stability under small deformations.

       The examples we provide are all compact quotients of nilpotent
       Lie groups.  To show non-existence, we use the "symmetrization"
       process and reduce the problem to non-existence of invariant
       structures with the given property. Then previous works by
       Dotti, Fino, Salamon and Parton \cite{DF1, FPS}
       lead to the results.

    The paper is organized as follows: In Section 2 we provide the
    material about the symmetrization of geometric structures on
    compact quotients of Lie groups. In Section 3 we provide a
    2-step nilmanifold example of hypercomplex manifolds which do not admit HKT
    structure and a one parameter family of hypercomplex structures,
     such that for one value of the parameter it
    admits HKT structure, while for the other values it does not.
    Finally, in Section 4,  we first prove that a compact Hermitian
    manifold with holomorphically trivial canonical bundle and
    vanishing Ricci tensor of its Bismut connection is (globally) conformally
    balanced. Then, we use the results from Section 2 to show that
    in a particular 2-parameter family of complex structures on the
Iwasawa manifold
    for all values of the parameter except one it does not admit
balanced metrics.

\section{Invariant structures on compact quotients of Lie groups}

In this section we prove some general facts about "symmetrizing"
geometric structures on compact quotients of a Lie group. We
start with the following theorem, which is a consequence of
\cite[Lemma 6.2]{Mil} where it is shown that {\it any
simply-connected Lie group which admits a discrete subgroup with
compact quotient is unimodular and in particular admits a
bi-invariant volume form $d\mu$}. For such manifolds we have:
\begin{teo}
Suppose that $M=\Gamma\backslash G$ is a compact quotient of a
simply-connected Lie group with a Riemannian metric $g$. Denote by
$d\mu$ a bi-invariant volume form and define a new
(left-invariant) metric by
$$\overline{g}(A,B)=\int_M g_m(A_m,B_m)d\mu. $$ Then we have the
following:

(a) $\overline{g}(\overline{\nabla}_A B, C) = \int_{m\in M
}g_m(\nabla_A B|_m,C_m)d\mu,$
where $A,B,C$ are projections of left-invariant vector fields from
$G$ to $M$.

     Similarly if  $B_1,...,B_k$ are projections of left-invariant
vector fields from $G$ to $M$ and
$$\overline{\omega}(B_1,...,B_k) = \int_{m \in M} \omega_m
(B_1|_m,...,B_k|_m)d\mu$$ where $\omega$ is a k-form on $M$, then
we have

(b) $d\overline{\omega}(B_1,..., B_{k+1})= \int_{m \in M}
d\omega|_m(B_1|_m,...B_{k+1}|_m)d\mu.$
\end{teo}

\begin{proof}
We follow the lines of \cite{Bel}. First note that
$$ \int_M Xf d\mu = \int_M \mathcal{L}_X (fd\mu) =
\int_Md(\it{i}_X fd\mu) = 0,$$ for any function $f$ and left-invariant
vector field $X$. Next, we remind the formula for the Levi-Civita
connection:
$$
\begin{array} {l l l}
     g(\nabla_X Y,Z)&=&1/2(Yg(X,Z)-Zg(X,Y)+Xg(Y,Z)+ \\
     & &g( \l
     X,Y \r ,Z+g( \l Z,X \r ,Y)-g( \l Y,Z \r ,X)).
\end{array}
$$
Now, using the above identity for left-invariant fields $A,B,C$, we
have:
$$
\begin{array}{ccl}
\int_M g_m(\nabla_A B|_m,C_m) d\mu  &= & \int_M g( \l
     A,B \r_m ,C_m)+g( \l C,A \r_m ,B_m)-g( \l B,C \r_m ,A_m))d\mu \\
     & = & \overline{g}( \l
     A,B \r ,C)+\overline{g}( \l C,A \r ,B)-\overline{g}( \l B,C \r
     ,A)\\
     &=& \overline{g}(\overline{\nabla}_A B,C)
\end{array}
$$
which proves (a). Using the definition for the exterior
differential on k-forms, we obtain (b) in a similar fashion.\end{proof}

As a consequence, we have the following:

\begin{teo} Suppose that $M$ is as in Theorem 2.1 and admits a
left-invariant complex structure $J$. Suppose moreover that $F$ is
a K\"ahler form of a non-invariant Hermitian metric $g$. Then the
left-invariant form defined as:
$$
\alpha(A_1,...,A_{2n-2}) = \int_M
F^{n-1}|_m(A_1|_m,...A_{2n-2}|_m) d\mu,
$$
for left-invariant vector fields, is equal to $\overline{F}^{n-1}$
for some K\"ahler form of a left-invariant Hermitian
metric. Moreover if $dF^{n-1} =0$ then $d\overline{F}^{n-1}=0$.
\end{teo}

\begin{proof}

Since $\alpha$ is a strictly positive  $(n - 1, n - 1)$-form, by a
linear algebra argument (see \cite{Mic}, p.279) it can be written
as $\overline{F}^{n-1}$, where $\overline F$ is a strictly
positive $(1, 1)$-form. Then it is a K\"ahler form of an invariant
metric and the last assertion follows from Theorem 2.1.
\end{proof}

\section{Examples of hypercomplex manifolds which do not admit HKT metric}
    A hypercomplex structure on a manifold $M$
is a triple of complex structures $\{ J_i \}_{i = 1,2,3}$
satisfying the quaternions relations $J_i^2 = - I$, $i = 1,2,3$,
$J_1 J_2 = - J_2 J_1 = J_3$. Let $g$ be a Riemannian metric on
$M$. Then $M$ is said to be hyper-Hermitian if it is Hermitian
with respect to every $J_i$, $i=1,2,3$.

A hyper-Hermitian manifold $(M, \{ J_i \}_{i = 1,2,3})$ is an HKT
manifold
if there is a connection $\nabla$ such
that
$$
\nabla g = 0, \quad \nabla J_i = 0, i = 1, 2,3,  \quad c(X, Y, Z) = g
(X, T (Y, Z)) \, \hbox{a 3-form}.
$$
If such connection exists it is unique. A hyper-Hermitian manifold
$M$ will admit a HKT connection if and only if $J_1 dF_1 = J_2
dF_2 = J_3 dF_3$, where $F_i$ is the K\"ahler form associated to
$(J_i, g)$, \cite{GP}.

      A hypercomplex structure on a Lie algebra
$\mathfrak g$ is a triple of complex structures
$\{ J_i
\}_{i = 1,2,3}$ satisfying the quaternions relations $J_i^2 = - I$,
$i = 1,2,3$, $J_1 J_2 = - J_2 J_1 = J_3$.

The hypercomplex structure will be called abelian if
$$
[J_i X, J_i Y] = [X, Y], \quad X, Y \in \mathfrak g.
$$
Let $g$ a inner product compatible with the hypercomplex structure.

By \cite{DF1}  $(\{ J_i \}_{i = 1,2,3}, g)$ is an HKT structure on
$\mathfrak g$ if and only if $g$ satisfy the extra
condition
$$
\begin{array} {l}
g ([J_1 X, J_1 Y], Z) + g ([J_1 Y, J_1 Z], X) + g ([J_1 Z, J_1 X], Y)  =\\
g ([J_2 X, J_2 Y], Z) + g ([J_2 Y, J_2 Z], X) + g ([J_2 Z, J_2 X], Y)  =\\
g ([J_3 X, J_3 Y], Z) + g ([J_3 Y, J_3 Z], X) + g ([J_3 Z, J_3 X], Y),
\end{array}
$$
for any $X, Y, Z \in {\mathfrak g}$.

When the hypercomplex structure is abelian the previous condition is
always satisfied for any  inner product $g$ compatible
with the hypercomplex structure.

We say that  a Lie group $G$ has a invariant HKT structure $( \{
J_i \}_{i = 1,2,3}, g)$ if the hypercomplex structure $\{  J_i
\}_{i = 1,2,3}$ and the metric  $g$ arise from corresponding
left-invariant tensors.  In \cite{DF1} it is shown that, if $G$ is
a 2-step nilpotent Lie group, the hypercomplex structure of an
invariant HKT structure is abelian. Then one has the following:

\begin{lemm} A 2-step nilpotent Lie group $G$ with a non abelian
hypercomplex structure admits no  invariant HKT  metric
compatible with such hypercomplex structure.
\end{lemm}

If we restrict to the 8-dimensional case  all  simply connected
nilpotent Lie groups which carry invariant abelian hypercomplex
structures are described in \cite{DF2}. There are three such
groups and they are central extensions of the real, complex or
quaternionic Heisenberg group respectively:
$$
N_1 = \mathbb{R}^3 \times H_5, \quad N_2 = \mathbb{R}^2 \times
H_3^{\C}, \quad N_3 = \mathbb{R} \times H_7.
$$
As remarked in \cite{DF2} $N_1$ and $N_2$ can only carry abelian
hypercomplex structures, but $N_3$ has
a non abelian hypercomplex
structure so it admits no invariant HKT  metric.

The Lie algebra ${\mathfrak n}_3$ of $N_3$ has structure equations

$$
\begin{array} {l}
\; [e_1, e_2] = [e_3, e_4] = - e_6,\\
\; [e_1, e_3] = - [e_2, e_4] = - e_7,\\
\; [e_1, e_4] = [e_2, e_3] = - e_8
\end{array}
$$

and a non abelian hypercomplex structure is given by
$$
\begin{array}{l}
J_1 (e_1) = e_2, J_1 (e_3) = e_4, J_1 (e_5) = e_6, J_1 (e_7) = e_8,\\
J_2 (e_1) = e_3, J_2 (e_2) =  - e_4, J_2 (e_5) =   e_7, J_2 (e_6)
= - e_8.
\end{array}
$$
In \cite{DF3} it is shown that a 8-dimensional nilpotent Lie algebra
$\mathfrak g$
which admit a hypercomplex structure is 2-step, its first Betti
number satisfies $b_1 ({\mathfrak
g})\geq 4$ and there exists a decomposition ${\mathfrak g} =
{\mathfrak g}_1 \oplus {\mathfrak g}_2$ such that\newline
i)  $\dim
{\mathfrak g}_i = 4$;\newline
ii) ${\mathfrak g}_i$ is invariant with respect to the hypercomplex
structure;\newline
iii) $[{\mathfrak g},
{\mathfrak g}] \subset {\mathfrak g}_2 \subset {\mathfrak z}$,
(${\mathfrak z}$: center of $\mathfrak g$).

If $b_1 ({\mathfrak g}) = 4, 5$ the full description of all
8-dimensional nilpotent Lie algebras with a hypercomplex structure
is given in \cite{DF3}.  If $b_1 ({\mathfrak g}) = 6, 7$ the
hypercomplex structure is abelian.

      If $G$ is a simply-connected nilpotent Lie
group and, if the structure equations of its Lie algebra are
rational, then there exists a discrete subgroup $\Gamma$ of $G$
for which $M =
\Gamma \backslash G$ is compact
\cite{Mal}. Any invariant HKT structure on $G$ will pass to a HKT
structure on $M$.

\begin{teo} Any compact quotient $ M= \Gamma \backslash G$ of a
2-step nilpotent Lie group $G$ with a non abelian
left-invariant hypercomplex structure
$\{ J_i \}_{i = 1, 2, 3}$ admits no HKT metric compatible with such
hypercomplex structure.
\end{teo}

\begin{proof}
Suppose that $M$ admits a  (non invariant) HKT metric $g$ compatible
with the non abelian left-invariant hypercomplex structure
$\{ J_i
\}_{i = 1, 2, 3}$. We will construct now, by integration, a  HKT
metric on $G$ which is left-invariant and so obtain a
contradiction with the previous Lemma.

$G$ admits a bi-invariant volume element $d \mu$, which induces one
on the compact quotient $M$ and possibly after rescaling,
we can suppose that the volume of $M$ is equal to 1.

Let $\Omega_i$ the K\" ahler form associated to $(J_i, g)$, then the
condition that  $g$ is a HKT metric is equivalent to
$$
J_1 d \Omega_1 = J_2 d \Omega_2 = J_3 d \Omega_3.
$$
Define
$$
F_i (A, B) := \int_{m \in M} (\Omega_i)_m (A_m, B_m) d \mu, \,  i = 1,
2, 3,$$ where $A, B$ are projections on $M$ of left-invariant
vectors fields on $G$. Then, since $J_i$ is left-invariant it is
easy to check that  $F_i$ is $J_i$-invariant, i.e. the K\" ahler
form of a positive definite metric $g_0$ and that it is
left-invariant and  defined as in Theorem 2.1.

    From Theorem 2.1 we conclude that
$$
d F_i (A, B, C) = \int_M d \Omega_i (A, B, C) d \mu.
$$
Then one has
$$
\begin{array}{c}
J_1 d F_1 (A, B, C) = -dF_1 (J_1 A, J_1 B, J_1 C) =  -\int_M d \Omega_1 (J_1 A,
J_1 B, J_1 C) d \mu\\
     = \int_M J_1 d \Omega_1 (A, B, C) d \mu = \int_M J_2 d \Omega_2 (A,
B, C) d \mu =
-\int_M d \Omega_2 (J_2 A, J_2 B, J_2 C) d \mu\\
     = - d F_2 (J_2 A, J_2 B, J_2 C) = J_2 d F_2 (A, B, C),
\end{array}$$
so  $g_0$ is a left-invariant HKT metric on $G$, which is impossible.
\end{proof}

Isomorphism classes of hypercomplex structures on 8-dimensional
nilpotent Lie groups are studied in \cite{La}. We show that there
exists on a 2-step nilpotent Lie group a one parameter family of
hypercomplex structures  which is abelian for one value of the
parameter and non-abelian for the others. Consider the family of
Lie algebras ${\mathfrak g}_t$ with structure equations:

$$
\begin{array}{l}
\; [e_1, e_2] = - t e_6, [e_3, e_4] = (1 - t) e_6, \\
   \; [e_1, e_3] = - t e_7, [e_2, e_4] = (t - 1) e_7,\\
\; [e_1, e_4] = - t e_8, [e_2, e_3] = (1 - t) e_8
\end{array}
$$
and the hypercomplex structure
$$
\begin{array} {l}
J_1 e_1 = e_2, J_1 e_3 = e_4, J_1 e_5 = e_6, J_1 e_7 = e_8,\\
J_2 e_1 = e_3, J_2 e_2 = - e_4, J_2 e_5 = e_7, J_2 e_6 = - e_8.
\end{array}
$$
For any $t \neq 0, 1$, ${\mathfrak g}_t$ is isomorphic to the Lie
algebra ${\mathfrak n}_3$ of  $N_3$. Indeed, one has that
${\mathfrak g}_t = {\mathfrak g}_1 \oplus {\mathfrak g}_2$, with
${\mathfrak g}_2 = {\mathfrak z} = \mbox {span} \{ e_5, e_6, e_7,
e_8 \}$ the center of ${\mathfrak g}$. Then if one considers the
endomorphism $J_Z$ ($Z \in {\mathfrak g}_2$) of ${\mathfrak g}_1$,
defined  by
\begin{equation}
   g (J_Z X, Y) = g ([X, Y], Z),
\end{equation}
where $g$ is the  inner product on ${\mathfrak g}_t$ such that the
basis $\{ e_1, \ldots, e_8 \}$ is orthonormal, one can check that
$J_Z$ is invertible for any non zero $Z \in {\mathfrak g}_2$
(compare also \cite[Example 7.4]{La}).

The hypercomplex structure
is abelian for $t = 1/2$ and
  non-abelian otherwise. Since the Lie
algebras are isomorphic for all  $t \neq 0, 1$, we obtain a
one-parameter family of invariant
hypercomplex structures on the compact quotient $M$ of $N_3$.
Moreover they are not equivalent since one is abelian and the
others are non-abelian. This proves the following:
\begin{teo}
The HKT structure is not stable under deformations, i.e. there
exists a small hypercomplex deformation of HKT structure which is
not HKT.
\end{teo}

\begin{remark} M. Verbitsky \cite{V1} informed us that stability of
HKT structures
under deformations holds for a large class of hypercomplex
manifolds which include $S^1\times S^{4n-1}$ and $SU(3)$. He
showed that any compact HKT manifold  for which $\Omega =
F_2+iF_3$ is (2,0)-part (with respect to $J_1$) of a closed 2-form
is HKT stable under small deformations. This follows from his
characterization of HKT manifolds as hypercomplex manifolds with
strictly q-positive $\partial$-closed (2,0)-form \cite{V2}. In
fact strictly q-positivity is an open condition and if
$\Omega=P^{(2,0)}\alpha$ where $P^{(2,0)}$ is the projection onto
(2,0)-forms and $d\alpha=0$, then $\Omega_t=P_t^{(2,0)}\alpha$
provides strictly q-positive $\partial_t$-closed form for deformed
hypercomplex structure $J_1^t,J_2^t,J_3^t$ for small $t$.

\end{remark}

\section{Compact complex nilmanifolds which do not admit a metric
with vanishing Ricci tensor of the Bismut connection}

    Let $(M,g,J)$ be a Hermitian manifold with K\"ahler form $F$.
    Then the Bismut connection is the only connection with torsion
    $T^B$ given by:
    $$ g(X,T^B(Y,Z)) = J dF(X,Y,Z) = -dF(JX,JY,JZ). $$
There is also a Chern connection with torsion $T^C$ uniquely
determined by
$$ g(X,T^C(Y,Z)) = dF(JX,Y,Z). $$ Both connections have Ricci
forms, representing up to a constant the first Chern class of the
complex structure. Then, as a consequence of formula (2.7.6) in
\cite{Ga1}, the relation between the two Ricci forms is:
\begin{equation}
Ric^B = Ric^C + d\delta F,
\end{equation}
where $\delta F$ is the co-differential (the adjoint of the
differential) of $F$. Moreover from (24) in \cite{Ga2} we have the
formula for the conformal change $\widetilde{g}=e^fg$:
$$
d\widetilde{\delta}\widetilde{F} = d\delta F+ dJd f,
$$
where by definition $Jdf(X) = -df(JX)$.  For $d\delta F$ we also
have:

\begin{lemm}
For a  compact Hermitian manifold $(M, J, g)$ with K\"ahler form
$F$ with co-differential $\delta F$ , $d\delta F =0$ iff M  is
balanced (i.e. $\delta F = 0$ or equivalently $dF^{n-1}=0$).
\end{lemm}

\begin{proof}

If $<,>$ is the extension of the metric $g$ on 2-forms, then
$$ 0 = \int_M <d\delta F,F> vol = \int_M <\delta F,\delta F> vol $$
and the result follows.

\end{proof}

    From here there is the following:
\begin{teo}
Suppose that a compact complex manifold $(M, J)$  admits a
holomorphic non-vanishing $(n,0)$-form. Then, if the Ricci tensor
of the Bismut connection of some Hermitian metric $g$ vanishes,
$(M,J,g)$ is conformally balanced and in particular admits a
balanced metric.
\end{teo}

\begin{proof}
First notice that, for the Ricci curvature of the Chern connection
on $M$, we have the well known formula $$ Ric^C = dJd \log det
(g_{\alpha \overline{\beta} }). $$ Let $\sigma$ denotes the
holomorphic (n,0)-form. It is locally represented as $$\sigma =
fdz^1\wedge \ldots\wedge  dz^n,$$ for a holomorphic function $f$.
Then $|\sigma|^2=|f|^2det(g_{\alpha\overline{\beta}})$ and
$$dJd\log |f|^2det(g_{\alpha\overline{\beta}})= dJd\log
det(g_{\alpha\overline{\beta}})$$ so
$$Ric^C=-dJd\log |\sigma|^2. $$

     Now from (2) above it follows
     that  $d\delta F = dJd log |\sigma|^2$ and after a
     conformal change $F\rightarrow |\sigma|^2F$ we can make $d \delta
F$ equal to 0.

Then from the Lemma above the new metric is balanced.
\end{proof}

Note that on any $2n$-dimensional compact quotient  of a nilpotent
Lie group  with an invariant complex structure there exists an
invariant closed $(n,0)$-form. For the next examples we need to
use the definition of Lie algebra in terms of structure equations
for the differentials of invariant 1-forms. It is equivalent to
the standard one giving the commutators of vector fields.
Similarly we can define complex structures as endomorphisms of the
space of invariant 1-forms. Consider now the family of Lie
algebras $ {\mathfrak g}_{s,t}$ defined by:

$$
\begin{array}{l}
d e ^i  = 0, i = 1, \ldots, 4\\
d e^5 =  s e^1 \wedge e^2 + 2s  e^3 \wedge e^4 + t e^1\wedge e^3 -t
e^2\wedge e^4,\\
d e^6 =  t e^1 \wedge e^4 + t e^2 \wedge e^3
\end{array}
$$
and the complex structure $J$
$$
Je^1 = e^2, Je^3 = e^4, Je^5 = e^6.
$$
For the corresponding $(1, 0)$-forms
$$
\omega_1 = e^1 + i e^2, \omega_2 = e^3 + i e^4, \omega_3 = e^5 + i e^6
$$
one has that
     $$
\begin{array} {l}
d \omega_i = 0, i = 1, 2,\\
d \omega_3 = \frac{1}{2} i s \, \omega_1 \wedge \overline \omega_1 +
i s \, \omega_2 \wedge \overline
\omega_2  + t \,  \omega_1 \wedge \omega_2
\end{array}
$$
and thus
$d(\omega_1\wedge \omega_2 \wedge \omega_3) = 0$, so there is
     a holomorphic (3,0)-form on the corresponding simply connected
     Lie group.

     For any $t\neq 0$, the Lie
     algebras ${\mathfrak g}_{s,t}$ are isomorphic to the real Lie
algebra underlying  the complex Heisenberg group
$$
G = \left\{ \left( \begin{array} {ccc} 1& z^1&z^3 \\ 0&1&z^2\\ 0&0&1
\end{array}
\right) : z^i \in \C \right\}
$$
with structure
equations
$$
\begin{array} {l}
d e^i = 0, i = 1, \ldots, 4,\\
d e^5 = e^1 \wedge e^3 - e^2 \wedge e^4,\\
d e^6 = e^1 \wedge e^4 + e^2 \wedge e^3,
\end{array}
$$
since ${\mathfrak g}_{s, t} = {\mathfrak g}_1 \oplus {\mathfrak g}_2$ with
$$
\begin{array} {l}
{\mathfrak g}_1 \cong \R^4 = {\mbox {span}} \{ e_1, e_2, e_3, e_4\},\\
{\mathfrak g}_2 \cong \R^2 = {\mbox {span}}  \{ e_5,  e_6 \}
\end{array}
$$
and  the endomorphism $J_Z$ defined by (1) of ${\mathfrak g}_1$ is
non-singular, for any non zero $Z \in {\mathfrak g}_2$ (compare
also \cite[Example 6.1 and 6.5]{La}).

Now for any invariant Hermitian metric its K\"ahler form is
of the form:
$$
\begin{array} {ll}
F= &x_1 i\, \omega_1  \wedge \overline \omega_1 + (x_2 + i x_3) \,
\omega_1  \wedge \overline \omega_2 + (x_2 - i x_3) \,
\overline \omega_1 \wedge \omega_2\\
{}& + \,  (x_4 + i x_5) \omega_1  \wedge \overline \omega_3 + \,
(x_4 - i x_5 )
\overline \omega_1 \wedge \omega_3 + x_6 i \, \omega_2 \wedge
\overline \omega_2\\
{}&+ (x_7 + i x_8) \, \omega_2 \wedge \overline \omega_3 + (x_7 - i
x_8) \, \overline \omega_2 \wedge \omega_3 + x_9 i
\, \omega_3
\wedge \overline \omega_3,
\end{array}
$$
where $x_i, i = 1, \ldots, 9$, are real numbers such that the
restriction of $F$ to any complex line is non zero. In
particular one has that $x_1, x_6, x_9 > 0$.
The invariant Hermitian metric is balanced if and only if $F$ is
orthogonal to the image of $d$ in $\Lambda^2 {\mathfrak g}^*$
(see \cite{AGS,FPS}).
The condition that $F$ is orthogonal to $d \omega_3$ becomes $s (x_1
+ \frac{1}{2} x_6) = 0$.
     Then, if $s \neq 0$,  there is no invariant
balanced Hermitian metric on ${\mathfrak g}_{s,t}$. If $s =0$ the
complex structure is the bi-invariant complex structure
$J_0$ on the complex Heisenberg group.  In this way the Iwasawa manifold
$M = \Gamma \backslash G$ (compact quotient of the complex Heisenberg
group by the discrete subgroup $\Gamma$
for which $z^i$
are Gaussian integers) comes equipped with a family
$J_{s,t}$ of invariant complex structures, for
$t\neq 0$, $s,t$ real and we have:
\begin{teo}
The Iwasawa manifold $(M,J_{s,t})$ admits a metric with vanishing Ricci
tensor of the Bismut connection if and only if $s=0$.
\end{teo}

\begin{proof}
Suppose that $s \neq 0$. If there exists a metric with vanishing
Ricci tensor of the Bismut connection, then by Theorem 4.1 there
exists a balanced metric, since there is a holomorphic (3,0)-form.
Now, by Theorem 2.2, there exists an invariant balanced metric.
However this is impossible as we showed earlier. Now when $s=0$,
there is an
  invariant balanced metric $g$. For any such metric
$Ric^C=Ric^B = dd^c log|\sigma|^2$ as in the proof of Theorem 4.1.
Then the conformal change $\tilde{g} = e^{|\sigma|}g$ (up to a
constant) provides a metric with vanishing Ricci tensor of the
Bismut connection. The  result also follows from Proposition 6.1
in \cite{FPS}.
\end{proof}

\begin{remark} With the above result we provide a counter-example of the
Conjecture 1.1 in \cite{GIP}. Moreover it shows that the property
"vanishing Ricci tensor of the Bismut connection" is not stable
under small deformations.
\end{remark}

\begin{remark}
Similar deformation appears in \cite{AB} where the authors prove
that the existence of balanced metric is not stable under small
deformations. The arguments presented here however are different.
\end{remark}

{\bf Acknowledgements} We are grateful to G.Papadopoulos and
M.Verbitsky for their remarks. G.G. wishes to thank "Abdus Salam"
ICTP, Trieste, where the work on this paper started in the Summer
of 2002.

\bigbreak 
\renewcommand{\thebibliography}{\list{\arabic{enumi}.\hfil}
{\settowidth\labelwidth{18pt}\leftmargin\labelwidth\advance
\leftmargin\labelsep\usecounter{enumi}}\def\newblock{\hskip.05em}
\sloppy\sfcode`\.=1000\relax}\newcommand{\bi}{\vspace{-3pt}\bibitem}

Authors address: G. Grantcharov, Department of Mathematics, Florida
International University, U.S.A. E-mail: grantchg@fiu.edu.

A. Fino:  Dipartimento di Matematica, Universit\`a  di
Torino, Via Carlo Alberto 10, Torino, Italy. E-mail: annamaria.fino@unito.it.

\end{document}